\documentclass[11pt]{article} 

\usepackage{amsmath}
\usepackage{amssymb}
\usepackage{amsthm}
\usepackage{amsfonts}
\newtheorem{Theorem}{\indent Theorem}[section]
\newtheorem{Proposition}{\indent Proposition}[section]
\newtheorem{Definition}[Theorem]{\indent Definition}
\newtheorem{Corollary}[Theorem]{\indent Corollary}
\newtheorem{Lemma}[Theorem]{\indent Lemma}
\theoremstyle{definition}
\newtheorem{example}[Theorem]{Example}

\def \qed{\hfill{\hbox{$\square$}}}
\def \proof{{\it Proof.  }}
\def \eproof{\qed\vskip .3cm}

\numberwithin{equation}{section}

\begin{document}
\title{Pseudo--spherical submanifolds with 1--type pseudo--spherical Gauss map}
\author{Burcu Bekta\c{s}, Elif \"{O}zkara Canfes and U\u{g}ur Dursun}

\date{}

\maketitle

\begin{abstract}
In this work, we study the pseudo--Riemannian submanifolds of a pseudo--sphere
with 1--type pseudo--spherical Gauss map. 
First, we classify the Lorentzian surfaces in a 4--dimensional 
pseudo--sphere $\mathbb{S}^4_s(1)$ with index $s$, $s=1, 2$, and having harmonic 
pseudo--spherical Gauss map. 
Then we give a characterization theorem for pseudo--Riemannian submanifolds of a pseudo--sphere
$\mathbb{S}^{m-1}_s(1)\subset\mathbb{E}^m_s$ with 1--type pseudo--spherical Gauss map,
and we classify spacelike surfaces and Lorentzian surfaces in the de Sitter space
$\mathbb{S}^4_1(1)\subset\mathbb{E}^5_1$ with 1--type pseudo--spherical Gauss map. 
Finally, according to the causal character of the mean curvature vector 
we obtain the classification of submanifolds of a pseudo--sphere 
having 1--type pseudo--spherical Gauss map with nonzero constant component in its spectral decomposition
\end{abstract}



\section{Introduction}
The theory of finite type submanifolds of a Euclidean space  
was initiated by B.-Y. Chen in late seventies and it has become
useful tool for investigation of submanifolds such as 
minimal submanifolds and parallel submanifolds.  
Since that time there have been rapid developments of finite type submanifolds of a Euclidean 
and a pseudo--Euclidean space, \cite{C3, C1, C2, Yoon}. 
Then, B.-Y. Chen and P. Piccinni extended the notion of finite type submanifolds to
differentiable maps on submanifolds, in particular to Gauss map, \cite{CP}.

A smooth map $\phi:M\longrightarrow\mathbb{E}^m$ of a compact
Riemannian manifold $M$ into a Euclidean space $\mathbb{E}^m$
is said to be of finite type if it has a finite spectral decomposition:
\begin{equation}
\label{specdec}
\phi=\phi_0+\phi_1+\phi_2+\cdots+\phi_k,
\end{equation}
where $\phi_0$ is a constant vector in $\mathbb{E}^m$,
and $\phi_i$'s are nonconstant $\mathbb{E}^m$--valued maps such that
$\Delta\phi_i=\lambda_{p_i}\phi_i$ with $\lambda_{p_i}\in\mathbb{R}$,
$i=1,\dots,k$. 
When the spectral resolution contains exactly $k$ nonconstant terms,
the map $\phi$ is called of k--type. Otherwise, it is of infinite type. 
Since $M$ is compact, $\phi_0$ is the center of mass of $M$ in $\mathbb{E}^m$.
On the other hand, one cannot make the spectral decomposition of a map 
on a noncompact Riemannian manifold in general. 
However, the notion of finite type immersions on noncompact manifolds was given as above
in \cite[page 124]{C2}. 
For the noncompact case the Laplacian operator may have zero eigenvalue for a nonconstant 
map in the spectral decomposition. 

Let $\textbf{x} : M \rightarrow \mathbb E^m_s$  be an isometric immersion 
from an $n$--dimensional pseudo--Riemannian manifold $M$ into a pseudo--Euclidean space $\mathbb E^m_s$. 
Let $G(n, m)$ denote the Grassmannian manifold consisting of all oriented $n$--planes of  $\mathbb E^m_s$. 
The {\it classical Gauss map} $\nu:M \rightarrow G(n, m)$ associated with 
$\textbf{x} $ is the map which carries each point $p \in M$ to the oriented $n$--plane of $\mathbb E^m_s$ 
obtained by parallel displacement of the tangent space $T_p M$ to the origin of $\mathbb E^m_s$. 
Since $G(n,m)$ can be canonically imbedded in the vector space 
$\bigwedge^n \mathbb E^m_s\cong\mathbb E^N_q$ for some integer $q$,
the classical Gauss map $\nu$ gives rise to a well--defined map from $M$ 
into the pseudo--Euclidean space $\mathbb E^N_q$, 
where $N= {m\choose n}$ and  $\bigwedge^n \mathbb E^m_s$
is the vector space obtained by the exterior products of $n$ vectors in $\mathbb E^m_s$, \cite{Kim}.

An isometric immersion from an $n$--dimensional Riemannian manifold $M$ into a Euclidean $(m-1)$--sphere 
$\mathbb  S^{m-1} $ can be viewed as one into a Euclidean $m$--space, and therefore 
the Gauss map associated with such an immersion can be determined in the ordinary sense. 
However, for the Gauss map to reflect the properties 
of the immersion into a sphere, instead of  into the Euclidean space, 
M. Obata modified the definition of Gauss map appropriately, \cite{Obata}.

Let $\textbf{x}: M \rightarrow \widetilde M$  be an isometric immersion 
from an $n$--dimensional Riemannian manifold $M$ 
into an $m$-- dimensional simply connected complete space $\widetilde M$ of constant curvature.
The generalized Gauss map in the Obata's sense is a map which assigns to each $p\in M$ the totally geodesic 
$n$--space tangent to $\textbf{x}(M)$ at $\textbf{x}(p)$.
In the case $\widetilde M^m = \mathbb  S^{m}$ (or resp. $\widetilde M^m = \mathbb  H^{m}$)
the generalized Gauss map is also called the {\it spherical Gauss map} 
(or resp. the {\it hyperbolic Gauss map}).
For spherical submanifolds, the concept of spherical Gauss map is more relevant than the classical one,
and also the geometric behavior of classical and spherical Gauss map are different. 
For example, the classical Gauss map of every compact Euclidean submanifold is mass--symmetric;
but the spherical Gauss map of a spherical compact submanifold is not mass--symmetric in general.

Later, in \cite{Ishihara} T. Ishihara studied Gauss map in a generalized sense of pseudo--Riemannian 
submanifolds of pseudo--Riemannian manifolds that also gives the Gauss map in the Obata's sense.

Let $\widetilde M^{m-1}_s$ denote the pseudo--sphere $\mathbb  S^{m-1}_s (1) \subset \mathbb E^m_s $ 
or the pseudo--hyperbolic space $\mathbb  H^{m-1}_s (-1) \subset \mathbb E^m_{s+1}$.
Let $\textbf{x}: M^n_t \rightarrow \widetilde M^{m-1}_s$  be an oriented 
isometric immersion from an $n$--dimensional pseudo--Riemannian 
manifold $M^n_t$ with index $t$ into the $(m-1)$--dimensional 
complete pseudo--Riemannian space $\widetilde M^{m-1}_s$ of constant curvature.
The generalized Gauss map in the Obata's sense is a map associated to $\textbf{x}$
which assigns to each $p\in M^n_t$ a totally geodesic 
$n$--subspace of $\widetilde M^{m-1}_s$  tangent to $\textbf{x}(M^n_t)$ at $\textbf{x}(p)$.
Since the totally geodesic $n$--subspace of $\widetilde M^{m-1}_s$ tangent to 
$\textbf{x}(M^n_t)$ at $\textbf{x}(p)$ is the pseudo--sphere $\mathbb  S^{n}_t (1)$ 
(or resp. the pseudo--hyperbolic space $\mathbb  H^{n}_t(-1)$), 
it determines a unique oriented $(n+1)$--plane containing $\mathbb  S^{n}_t (1)$ 
(or resp. $\mathbb  H^{n}_t(-1)$). 
Thus, the generalized Gauss map in the Obata's sense can be extended to a map $\hat \nu$ of $M^n_t$ into 
the Grassmannian manifold $G(n+1, m)$ in the natural way, and the composition $\tilde \nu$  of $\hat \nu$ 
followed by the natural inclusion of $G(n+1, m)$ into a pseudo--Euclidean space 
$\mathbb E^N_q$, $N= {m\choose n+1}$, for some integer $q$ is the \textit{pseudo--spherical Gauss map} 
(or resp. the \textit{pseudo--hyperbolic Gauss map}).

In \cite{CP}, B.-Y. Chen and P. Piccinni obtained characterization 
and classification of submanifolds with 
1--type Gauss map, and they also gave the complete classification of minimal surfaces
of $\mathbb{S}^{m-1}$ with 2--type Gauss map. In \cite{Dursun}, 
U. Dursun studied hypersurfaces of hyperbolic space
with 1--type Gauss map, and he classified hypersurfaces of a hyperbolic space 
in Lorentz--Minkowski space $\mathbb{E}^m_1$ with at most two distinct
principal curvatures and 1--type Gauss map.

In \cite{CL}, B.-Y. Chen and H.-S. Lue studied spherical submanifolds 
with finite type spherical Gauss map, 
and they obtained some characterization and classification results. 
In \cite{Dur-Bek}, B. Bekta{\c s} and U. Dursun determined submanifolds of the unit sphere 
$\mathbb S^{m-1}$ with nonmass--symmetric 1--type spherical Gauss map, and  they also classified  
constant mean curvature surfaces in $\mathbb S^3$ with mass--symmetric 2--type spherical Gauss map.
Recently, the fundamental results about the submanifolds of hyperbolic space 
with finite type hyperbolic Gauss map were obtained in \cite{Dur-Yeg}.

In this paper, we study submanifolds of a pseudo--sphere with finite type 
pseudo--spherical Gauss map. First we determine the Lorentzian surfaces in a
pseudo--sphere $\mathbb{S}^4_s(1)\subset\mathbb{E}^5_s$ with harmonic
pseudo--spherical Gauss map. Then we completely classify spacelike surfaces 
and Lorentzian surfaces in the de Sitter space $\mathbb{S}^4_1(1)\subset\mathbb{E}^5_1$
with 1--type pseudo--spherical Gauss map. 
According to the casual character of the mean curvature vector we classify the submanifolds
of a pseudo--sphere with 1--type pseudo--spherical Gauss map having a nonzero
constant component in its spectral decomposition.
Also, we determine the marginally trapped spacelike surfaces in the de Sitter space 
$\mathbb{S}^4_1(1)\subset\mathbb{E}^5_1$ having 1--type pseudo--spherical Gauss map
with nonzero constant component in its spectral decomposition.  
Moreover, we prove that an $n$--dimensional pseudo horosphere in 
$\mathbb{S}^{n+1}_s(1)\subset\mathbb{E}^{n+2}_s$ 
has biharmonic pseudo--spherical Gauss map.

\section{Preliminaries}
Let $\mathbb{E}^m_s$ be the pseudo--Euclidean space equipped with the canonical pseudo--Euclidean metric
of index $s$ given by
\begin{equation}
\label{metric}
\tilde{g}=\sum_{i=1}^{m-s} dx_i^2-\sum_{j=m-s+1}^{m} dx_j^2,
\end{equation}
where $(x_1, x_2,\dots, x_m)$ is a rectangular coordinate system of $\mathbb{E}^m_s$.

A vector $v \in \mathbb{E}^m_s$ is called spacelike if $\langle v, v\rangle>0$ or $v=0$, 
timelike if $\langle v,v\rangle<0$, and  
lightlike(or null) if $\langle v, v\rangle=0$ and $v\neq 0$.

For $x_0\in\mathbb{E}^m_s$ and $c$ is a positive real number, we put 
\begin{eqnarray}
\nonumber
\label{pseudosphy}
\mathbb{S}^{m-1}_s(x_0,c)=\left\{x=(x_1,x_2,\dots, x_m)\in\mathbb{E}_s^{m}
\;|\; \langle x-x_0, x-x_0\rangle= c^{-1}\right\},\\
\nonumber
\mathbb{H}^{m-1}_{s-1}(x_0,-c)=\left\{x=(x_1,x_2,\dots, x_m)\in\mathbb{E}_{s}^{m}
\;|\; \langle x-x_0, x-x_0\rangle=-c^{-1}\right\},
\end{eqnarray}
where $\langle , \rangle$ denotes the indefinite inner product on $\mathbb{E}^m_s$.
Then, $\mathbb{S}^{m-1}_s(x_0,c)$ and $\mathbb{H}^{m-1}_{s-1}(x_0,-c)$ are complete
pseudo--Riemannian manifolds of constant curvature $c$ and $-c$ centered at $x_0$, 
called pseudo--sphere and pseudo--hyperbolic space, respectively. 
When $x_0$ is the origin and $c=1$ 
we simply denote $\mathbb{S}^{m-1}_s(x_0, 1)$ and $\mathbb{H}^{m-1}_{s-1}(x_0,-1)$ 
by $\mathbb{S}^{m-1}_s(1)$ and $\mathbb{H}^{m-1}_{s-1}(-1)$. 
The manifolds $\mathbb{E}^m_1$, $\mathbb{S}^{m-1}_1$ and $\mathbb{H}^{m-1}_1$ 
are known as the Minkowski, de Sitter, and anti--de Sitter spaces, respectively.

The light cone $\cal{LC}$ in $\mathbb{E}^m_s$ is defined by
\begin{equation}
\label{lightcone}
{\cal LC}=\{v\in\mathbb{E}^m_s\;|\; \langle v, v\rangle=0\}.
\end{equation}

Let $M$ be an oriented $n$--dimensional pseudo--Riemannian submanifold in a
pseudo--Euclidean space $\mathbb E^m_s$. 
We choose an oriented local orthonormal frame field $\{e_1, \dots, $ $ e_{m} \}$
on $M$ such that  $e_1,\dots,e_{n}$ are tangent to $M$, and $e_{n+1}, \dots, e_{m}$ are normal to $M$
with the signatures $\varepsilon_A=\langle e_A,e_A \rangle=\pm 1$, $A=1,\dots, m$.
We use the following convention on the range of indices unless mentioned otherwise:
$$1\leq A,B,C,\ldots\leq m, \; 1\leq i,j,k,\ldots \leq n, \; n+1\leq r,s, t,\ldots \leq m.$$

Let $\widetilde{\nabla}$ be the Levi--Civita connection of $\mathbb E^m_s$
and  $\nabla$  the induced connection on $M$. With respect to the frame field 
of $M$ chosen above, let $\{\omega_A\}$ be 
the dual frame  and $\{\omega_{AB}\}$ the connection forms 
associated to $\{e_A \}$ with $\omega_{AB}+\omega_{BA}=0$. 
Thus the Gauss and Weingarten formulas are given, respectively, by 
\begin{equation}
\nonumber
\widetilde{\nabla}_{e_k}e_i= \sum_{j=1}^{n}\varepsilon_j \omega_{ij}(e_k) e_j + \sum_{r=n+1}^{m} \varepsilon_r h^r_{ik} e_r 
\end{equation}
and
\begin{equation}
\nonumber
\widetilde{\nabla}_{e_k}e_r=  - A_r(e_k)+ \sum_{s=n+1}^{m} \varepsilon_s \omega_{rs}(e_k) e_s,
\end{equation}
where $h^r_{ij}$'s are the coefficients of
the second fundamental form $h$, and
$A_r$ is the Weingarten map in the direction $e_r$.

The mean curvature vector $H$, the squared length $\| h\|^2$ of
the second fundamental form $h$ and the scalar curvature $S$ of $M$ in $\mathbb{E}^m_s$ are defined,
respectively, by
\begin{align}
\label{meancurvector}
H &=\frac{1}{n}\sum_{r=n+1}^m\varepsilon_r\mbox{tr}A_r e_r,\\
\label{sqr-length}
\| h\|^2 &= \sum_{i,j=1}^n\sum_{r=n+1}^m \varepsilon_i\varepsilon_j \varepsilon_r  h^r_{ij}h^r_{ji},\\
\label{scalarcurv}
S &=n^2\langle H,H\rangle - \| h\|^2.
\end{align}

The Codazzi equation of $M$ in $\mathbb{E}^m_s$ is given by
\begin{align} 
\label{codazzi}
\begin{split}
&h^r_{ij,k}= h^r_{jk,i},\\
&h^r_{jk,i} = e_i(h^r_{jk}) -
\sum_{\ell=1}^{n} \varepsilon_{\ell}\left (  h^r_{\ell k} \omega_{j\ell}(e_i) +
h^r_{\ell j} \omega_{k \ell}(e_i) \right ) + \sum_{s=n+1}^{m} \varepsilon_s h^s_{jk} \omega_{sr} (e_i).
\end{split}
\end{align}
Also, from the Ricci equation of $M$ in $\mathbb{E}^{m}_s$, we have
\begin{equation} 
\label{normal-curv}
R^D(e_j, e_k; e_r, e_s) = \langle [A_r,A_s](e_j), e_k \rangle =
\sum_{i=1}^{n} \varepsilon_i \left (h^r_{ik} h^s_{ij} - h^r_{ij} h^s_{ik}\right),
\end{equation}
where $R^D$ is the normal curvature tensor associated to the normal connection
$D$ defined by $D_{e_k}e_r=\sum_{s=n+1}^{m} \varepsilon_s \omega_{rs}(e_k) e_s$.
A submanifold $M$ in $\mathbb E^m_s$ is said to have flat normal bundle 
if $R^D$ vanishes identically.

In particular, if $M$ is immersed in the pseudo--sphere $\mathbb{S}^{m-1}_s(c)\subset\mathbb{E}^m_s$, 
we have 
\begin{equation} 
\label{mean-and-sff}
H= \hat H - c{\bf x}, \quad h(X,Y) = \hat  h(X,Y)  -  c\left\langle X,Y\right\rangle {\bf x},
\end{equation}
where $\hat{H}$ and $\hat{h}$ are the mean curvature vector and 
the second fundamental form in $\mathbb{S}^{m-1}_s(c)$. 
Then, equation \eqref{scalarcurv} becomes 
\begin{equation} 
\label{scalar-curv-sphere}
S= cn(n-1) + n^2 |\hat H|^2-  \| \hat h\|^2.
\end{equation}

Following \cite{Dursun} we can have the definition of finite type maps
on pseudo--Riemannian submanifolds of pseudo--sphere or pseudo--hyperbolic space
as follows:
\begin{Definition}
A smooth map $\phi:M_t\rightarrow\mathbb{S}^{m-1}_s(c)\subset\mathbb{E}^m_s$
(or resp. $\phi:M_t\rightarrow\mathbb{H}^{m-1}_{s-1}(-c)\subset\mathbb{E}^m_s$) from 
a pseudo--Riemannian submanifold $M_t$ with index $t$ into a pseudo--sphere 
$\mathbb{S}^{m-1}_s(c)$ (or resp. into a pseudo--hyperbolic space $\mathbb{H}^{m-1}_{s-1}(-c)$)
is called of finite type in $\mathbb{S}^{m-1}_s(c)$ (resp. in $\mathbb{H}^{m-1}_{s-1}(-c)$) 
if the map $\phi$ has the following spectral decomposition:
\begin{equation}
\label{specde}
\phi=\phi_1+\phi_2+\cdots+\phi_k, 
\end{equation}
where $\phi_i$'s are nonconstant $\mathbb{E}^m_s$--valued maps on $M_t$
such that $\Delta\phi_i=\lambda_{p_i}\phi_i$ with 
$\lambda_{p_i}\in\mathbb{R}, \;\; i=1,2, \dots, k$.
If the spectral decomposition \eqref{specde} contains exactly $k$ nonconstant maps,
then the map $\phi$ is called of $k$--type, 
and when one of the $\lambda_{p_i}$'s vanishes, the map $\phi$ is called of null $k$--type.
\end{Definition}

Note that for a finite type map, one of the components in this spectral decomposition 
may still be constant. 

For a finite type map on compact submanifold, the minimal polynomial criterion
was given in \cite{CL}.  
However, if the manifold is not compact, then the existence of minimal polynomial 
does not imply that it is of finite type, \cite{CPe}. 
A criteria for finite type maps was given in \cite{Dur-Yeg} as follows:
\begin{Theorem}
\cite{Dur-Yeg}
Let $\phi:M_t\longrightarrow\mathbb{E}^m_s$ be a smooth map from a pseudo--Riemannian manifold $M_t$
with index $t$ into a pseudo-Euclidean space $\mathbb{E}^m_s$, 
and let $\tau=\mbox{div}(d\phi)$ be the tension field of $\phi$. Then,
\begin{itemize} 
\item[i.] If there is a nontrivial polynomial $Q$ such that $Q(\Delta)\tau=0$, 
then $\phi$ is either of infinite type or of finite type with type number $k\leq\mbox{deg}(Q)+1$;

\item[ii.] If there is a nontrivial polynomial $P$ with simple roots such that
$P(\Delta)\tau=0$, then $P$ is of finite type with type number $k\leq\mbox{deg}(P)$.
\end{itemize}
\end{Theorem}

\section{Pseudo--spherical Gauss map}
In this section, we give the details of the pseudo--spherical Gauss map
mentioned in Introduction and study the pseudo--Riemannian submanifolds 
with harmonic pseudo--spherical Gauss map.

Let $\mathbf{x}: M_t \rightarrow \mathbb{S}^{m-1}_s(1)\subset\mathbb{E}^{m}_s$
be an isometric immersion from an $n$--dimensional pseudo--Riemannian manifold 
$M_t$ with index $t$ into a pseudo--sphere $\mathbb{S}^{m-1}_{s}(1)\subset\mathbb{E}^{m}_{s}$.
The map $\hat{\nu}: M_t \rightarrow G(n+1, m)$, 
called pseudo--spherical Gauss map in Obata's sense, 
of an immersion $\mathbf{x}$ into the Grassmannian manifold  
$G(n+1, m)$ assigns to each point $p\in M_t$ the great pseudo $n$--subsphere
$\mathbb{S}^{n}_{t}(1)$ of $\mathbb{S}^{m-1}_{s}(1)$ 
tangent to $\mathbf{x}(M_t)$ at $\mathbf{x}(p)$.
The great pseudo $n$--subsphere $\mathbb{S}^{n}_{t}(1)$ of $\mathbb{S}^{m-1}_{s}(1)$
are naturally identified with the Grassmannian manifold of oriented $(n+1)$--planes
through the center of $\mathbb{S}^{m-1}_{s}(1)$ in $\mathbb{E}^{m}_{s}$ 
since such $(n+1)$--planes determine unique great pseudo $n$--subsphere and conversely.

On the other hand, since the Grassmannian manifold $G(n+1, m)$ can be canonically imbedded in a
pseudo--Euclidean space $\bigwedge^{n+1}\mathbb E^{m}_{s} \cong \mathbb E^{N}_q$ 
obtained by the exterior products of $n+1$ vectors in 
$\mathbb E^{m}_{s}$ for some positive integer $q$, 
the composite $\tilde{\nu}$ of $\hat{\nu}$ followed 
by the natural inclusion of $G(n+1,m)$ in $\mathbb{E}^{N}_{q}$ is also 
called the pseudo--spherical Gauss map, 
where $N= {m\choose n+1}$.

For each point $p \in M_t$, let  $ e_1,\dots, e_{n}$
be an orthonormal basis of $T_p M_t$ with the signatures 
$\varepsilon_i = \left\langle e_i, e_i\right\rangle = \pm 1, \; i=1, \dots, n$. 
Then, the $n+1$ vectors $\mathbf{x}(p), e_1,\dots, e_{n}$ determine a 
linear $(n+1)$--subspace in $\mathbb E^{m}_{s}$. 
The intersection of this linear subspace and $\mathbb S^{m-1}_s(1)$ is
a totally  geodesic pseudo $n$--sphere $\mathbb{S}^{n}_{t}(1)$  
determined by $T_p M_t$.

Let $\mathbb E^{m}_{s}$ be a pseudo--Euclidean space with
two orthonormal bases $\{f_1, f_2, \dots, f_{m}\}$ and
$\{g_1, g_2, \dots, g_{m}\}$. Let $f_{i_1} \wedge \cdots \wedge f_{i_{n+1}}$ and
$g_{j_1} \wedge \cdots \wedge g_{j_{n+1}}$ be two vectors 
in $\bigwedge^{n+1} \mathbb E^{m}_{s}$.
Define an indefinite inner product 
$\left\langle \left\langle, \right\rangle \right\rangle$ on
$\bigwedge^{n+1} \mathbb E^m_{s}$ by 
\begin{equation} \label{inner-prod}
\left\langle \left\langle f_{i_1} \wedge \cdots \wedge f_{i_{n+1}},  g_{j_1} \wedge \cdots \wedge g_{j_{n+1}}\right\rangle \right\rangle
= \det( \left\langle f_{i_\ell}, g_{j_k}\right\rangle).
\end{equation}
Therefore, we may identify $\bigwedge^{n+1} \mathbb E^{m}_{s}$ with 
some pseudo--Euclidean space $\mathbb E^N_q$ for some positive integer $q$,
where $N= {m\choose n+1}$.

For an immersion  $\mathbf{x}: M_t \rightarrow \mathbb S^{m-1}_s(1)\subset\mathbb E^{m}_{s}$, 
the pseudo--spherical Gauss map in the Obata's sense can be considered 
as $\hat \nu: M_t \rightarrow G(n+1, m)$
which carries each $p \in M_t$ to 
$\hat\nu(p) = (\mathbf{x} \wedge e_{1}\wedge \cdots \wedge e_{n}) (p) $.
Since $\left\langle  \left\langle \hat \nu, \hat \nu \right\rangle \right\rangle =  
\varepsilon_1 \varepsilon_2 \cdots \varepsilon_n = \pm 1$, 
$G(n+1, m)$ is a submanifold of 
$\mathbb S^{N-1}_q(1) \subset \mathbb E^N_q $ 
or \linebreak $\mathbb H^{N-1}_{q-1} (-1) \subset \mathbb E^N_q$. 
Thus, considering the natural inclusion of $G(n+1, m)$ into $\mathbb E^N_q$,
the pseudo--spherical Gauss map $\tilde \nu$ associated 
with $\mathbf{x}$ is   given by 
\begin{equation} \label{sph-gauss-map}
\tilde \nu = \mathbf{x} \wedge e_{1} \wedge e_{2} \wedge \cdots \wedge e_{n} : M_t  \rightarrow G(n+1,m) \subset \mathbb{E}^N_q.
\end{equation}

Now, by differentiating $\tilde\nu$ in \eqref{sph-gauss-map} we have
\begin{equation}
\label{dif1}
e_i\tilde\nu=\sum_{k=1}^n\sum_{r=n+1}^{m-1}\varepsilon_r h_{ik}^r {\bf{x}}\wedge e_1\wedge\cdots\wedge\underbrace{e_r}_{k-th}\wedge\cdots e_n.
\end{equation}
Using $\nabla_{e_i} e_i=\sum_{j=1}^n \varepsilon_j\omega_{ij}(e_i)e_j$ and \eqref{dif1}, we obtain
\begin{equation}
\label{dif2} 
\left(\nabla_{e_i}e_i\right)\tilde\nu=
\sum_{j,k=1}^n\sum_{r=n+1}^{m-1} \varepsilon_j\varepsilon_r h_{jk}^r\omega_{ij}(e_i) {\bf x}\wedge e_1\wedge\cdots\wedge\underbrace{e_r}_{k-th}\wedge\cdots\wedge e_n.
\end{equation}
Since the Laplacian of $\tilde\nu$ is defined by
\begin{equation}
\label{deflaplas}
\Delta\tilde\nu=\sum_{i=1}^n\varepsilon_i(\nabla_{e_i} e_i-e_ie_i)\tilde\nu,
\end{equation}
considering equations \eqref{dif1}, \eqref{dif2} and the Codazzi equation \eqref{codazzi} we obtain 
\begin{align}
\label{dif3}
\begin{split}
\Delta\tilde\nu=& \| h \|^2 \tilde\nu + n H\wedge e_1\wedge\cdots\wedge e_n \\
&-n\sum_{k=1}^n {\bf x}\wedge e_1\wedge\cdots\wedge\underbrace{D_{e_k}H}_{k-th}\wedge\cdots\wedge e_n\\
&-\sum_{i,j,k=1 \atop j\neq k}^n\sum_{r,s=n+1}^{m-1}\varepsilon_i\varepsilon_r\varepsilon_s h_{ik}^r h_{ij}^s {\bf x} \wedge e_1\wedge\cdots\wedge\underbrace{e_s}_{j-th}\wedge\cdots\wedge\underbrace{e_r}_{k-th}\wedge\cdots\wedge e_n,
\end{split}
\end{align}
where $D$ is the normal connection induced on $M_t$.

Therefore, by using \eqref{normal-curv} and \eqref{mean-and-sff} in \eqref{dif3} 
we have
\begin{Lemma}
Let $M_t$ be an n--dimensional oriented pseudo--Riemannian submanifold with index $t$
of a pseudo $(m-1)$--sphere 
$\mathbb{S}^{m-1}_s(1)\subset\mathbb{E}^m_s$. 
Then the Laplacian of the pseudo--spherical Gauss map 
$\tilde\nu: M_t\rightarrow G(n+1, m)\subset\mathbb{E}^N_q$, 
$N= {m\choose n+1}$, for some $q$ is given by
\begin{align}
\label{calc-laplce-lem}
\begin{split}
\Delta \tilde \nu = &  \| \hat h \|^2 \tilde \nu + n \hat H  \wedge e_{1} 
\wedge  \cdots \wedge e_{n} 
- n \sum_{k=1}^n   {\bf x} \wedge  e_{1} \wedge \cdots \wedge  
\underbrace{D_{e_k} \hat H }_{k-th} \wedge \cdots \wedge e_{n}  \\ 
& +\sum_{j,k=1 \atop j\neq k}^n\sum_{r,s=n+1,\atop r<s}^{m-1} \varepsilon_r \varepsilon_s R^r_{sjk} 
{\bf x} \wedge e_{1} \wedge \cdots 
\wedge  \underbrace{e_r}_{j-th}  \wedge \cdots \wedge  \underbrace{e_s}_{k-th} \wedge \cdots \wedge e_{n},
\end{split}
\end{align}
where $R^r_{sjk}=R^D(e_j, e_k; e_r, e_s)$.
\end{Lemma}

\begin{Proposition} 
\label{prop-1}
Let ${\bf x}:(M_t,g)\longrightarrow\mathbb{S}^{m-1}_s(1)\subset\mathbb{E}^m_s$ be an isometric immersion
from an $n$--dimensional pseudo-Riemannian manifold $M_t$ with the metric $g$ and index $t$ into a pseudo-sphere 
$\mathbb{S}^{m-1}_s(1)\subset\mathbb{E}^m_s$. Then we have the following statements:
\begin{itemize}
\item [i.] The pseudo-spherical Gauss map in the Obata's sense 
$\hat\nu:(M_t,g)\longrightarrow G(n+1,m)$ is a harmonic map
if and only if $M_t$ has zero mean curvature in $\mathbb{S}^{m-1}_s(1)$;

\item [ii.] The pseudo-spherical Gauss map $\tilde\nu:(M_t,g)\longrightarrow\mathbb{E}^N_q$ 
with $N= {m\choose n+1}$ and some positive integer $q$ 
is harmonic if and only if 
$M_t$ has zero mean curvature in $\mathbb{S}^{m-1}_s(1)$, flat normal connection 
and constant scalar curvature $S=n(n-1)$.
\end{itemize}
\end{Proposition}

The proof of i. is similar to the Proposition 3.2 given in \cite{CL}, 
and the proof of ii. comes from equation \eqref{calc-laplce-lem} directly. 
Also, it can be seen that $\|\hat{h}\|^2$ could be zero even if $\hat{h}\neq 0$. 
Thus, submanifold with harmonic pseudo--spherical Gauss map 
does not need to be totally geodesic in general. 

Now, we will discuss spacelike surfaces and Lorentzian surfaces in $\mathbb{S}^4_s$, $s=1,2$, with 
harmonic pseudo--spherical Gauss map. 

From  part ii of Proposition \ref{prop-1} we have 
\begin{Corollary}
A spacelike surface $M$ in $\mathbb{S}^4_2(1)\subset\mathbb{E}^5_2$ has harmonic pseudo-spherical Gauss map
if and only if $M$ is a totally geodesic surface in $\mathbb{S}^4_2(1)$.
\end{Corollary}

In \cite{C5}, B.-Y. Chen gave the classification of all minimal Lorentz surfaces of 
constant Gaussian curvature one 
in the pseudo--sphere with arbitrary dimension and index. 
In the given classification theorem there are 3 surfaces in which one of them has 
nonflat normal bundle. Therefore by using Proposition 5.1 in \cite{C5} and 
Proposition 3.1 we state
\begin{Theorem}
A Lorentzian surface $M_1$ in $\mathbb{S}^4_s(1)$, $s=1,2$, has harmonic pseudo spherical--Gauss map 
if and only if one of the following cases occurs:
\begin{itemize}
\item [i.] $M_1$ is an open portion of the totally geodesic de Sitter space
$\mathbb{S}^2_1(1)\subset\mathbb{S}^4_s(1)$, 

\item [ii.] The immersion $L:M_1\longrightarrow\mathbb{S}^4_s(1)\subset\mathbb{E}^5_s$ 
is given by
\begin{equation}
L(u,v)=\frac{z(u)}{u+v}-\frac{z'(u)}{2},
\end{equation}
where $z(u)$ is a spacelike curve with constant speed $2$ lying in the light cone $\cal{LC}$ satisfying 
$\langle z'', z''\rangle=0$ and $z'''\neq 0$.
\end{itemize}
\end{Theorem}

\section{Submanifolds with 1--type pseudo--spherical Gauss map}
In this section, we classify pseudo--spherical submanifolds in 
$\mathbb{S}^{m-1}_s(1)\subset\mathbb{E}^m_s$ 
with 1--type pseudo--spherical Gauss map.

\begin{Theorem}
\label{1typesubmanifold}
A pseudo--Riemannian submanifold $M_t$ with index t of \linebreak
$\mathbb{S}^{m-1}_s(1)\subset\mathbb{E}^m_s$ 
has 1--type pseudo--spherical Gauss map 
if and only if $M_t$ has zero mean curvature in $\mathbb{S}^{m-1}_s(1)\subset\mathbb{E}^m_s$, 
constant scalar curvature and flat normal connection. 
\end{Theorem}

\proof
If the pseudo--spherical submanifold of $\mathbb{S}^{m-1}_s(1)\subset\mathbb{E}^m_s$ 
has 1--type pseudo--spherical Gauss map $\tilde\nu$, 
then $\Delta\tilde\nu$ and $\tilde\nu$ are proportional. 
Thus, we see from \eqref{calc-laplce-lem} that $\tilde\nu$ is of 1--type 
if and only if $\hat{H}=R^r_{sjk}=0$ and $\|\hat{h}\|^2$ is constant. 
Hence, from \eqref{scalar-curv-sphere} it is seen that the scalar curvature is constant.
\eproof
  
The standard imbedding of $\mathbb{S}^1(2)\times\mathbb{S}^1(2)$ in $\mathbb{S}^3(1)$ is called the 
\textit{Clifford torus}. It is well known that the Clifford torus is flat and minimal in $\mathbb{S}^3(1)$.
 
\begin{Theorem}
\label{1typessurfacein41}
A non--totally geodesic spacelike surface $M$ in $\mathbb{S}^4_1(1)\subset\mathbb{E}^5_1$ 
has 1--type pseudo--spherical Gauss map if and only if 
it is an open portion of the Clifford torus lying fully in a totally geodesic 3--sphere 
$\mathbb{S}^3(1)\subset\mathbb{S}^4_1(1)$. 
\end{Theorem}

\proof
Suppose that a spacelike surface in $\mathbb{S}^4_1(1)\subset\mathbb{E}^5_1$ has 
1--type pseudo--spherical Gauss map. 
Following the result given in \cite{gorokh} and Theorem \ref{1typesubmanifold} it is seen that 
the only maximal surface of $\mathbb{S}^4_1(1)$ with constant Gaussian curvature and flat normal
connection is an open portion of the Clifford torus.

Conversely, suppose that $M$ is an open portion the Clifford torus 
in $\mathbb{S}^3(1)\subset\mathbb{S}^4_1(1)$.
Then it is easy to show that the pseudo--spherical Gauss map of the Clifford torus satisfies 
$\Delta\tilde\nu=2\tilde\nu$. Thus, it has 1--type pseudo--spherical Gauss map.
\eproof

The standard imbedding of $\mathbb{S}^1(2)\times\mathbb{S}^1_1(2)$ in $\mathbb{S}^3_1(1)$ 
is called the \textit{pseudo--Riemannian Clifford torus}.
We can imbed the pseudo--Riemannian Clifford torus in $\mathbb{S}^4_1(1)$ and give a parametrization
in $\mathbb{S}^4_1(1)$ as follows:
\begin{equation}
\label{prclifford}
r(u,v)=\frac{1}{\sqrt{2}}(0, \cos u, \sin u, \cosh v, \sinh v).
\end{equation}
By a direct calculation, we obtain that its pseudo--spherical Gauss map $\tilde\nu$ 
satisfies $\Delta\tilde\nu=2\tilde\nu$. 
Note that this surface has zero mean curvature in $\mathbb{S}^4_1(1)$, 
and the Gaussian and normal curvatures are zero.

Similarly, considering the result given in \cite{gorokh} and Theorem \ref{1typesubmanifold}
we can give the following classification 
theorem for Lorentzian surfaces in $\mathbb{S}^4_1(1)\subset\mathbb{E}^5_1$ 
with 1--type pseudo--spherical Gauss map.

\begin{Theorem}
\label{1typelsurfacein41}
A non--totally geodesic Lorentzian surface $M_1$ in $\mathbb{S}^4_1(1)\subset\mathbb{E}^5_1$ 
has 1--type pseudo--spherical Gauss map if and only if 
it is an open portion of the pseudo--Riemannian Clifford torus lying fully 
in a totally geodesic pseudo--sphere 
$\mathbb{S}^3_1(1)\subset\mathbb{S}^4_1(1)$.  
\end{Theorem}

Now, we discuss the pseudo--spherical submanifolds with 1--type pseudo--spherical
Gauss map such that its spectral decomposition contains a nonzero constant component.
For later use, we need the following lemma.

\begin{Lemma} 
\label{isothermal-surf-Gauss-map}
For a pseudo--Riemannian hypersurface $M_t$ with index $t$ of $\mathbb S^{n+1}_s(1)\subset \mathbb E^{n+2}_s$ 
we have
\begin{equation}
\label{laplce-lem-2}
\Delta ( e_{n+1}  \wedge e_{1} \wedge e_{2} \wedge \cdots \wedge e_{n}) = 
n \hat \alpha  \tilde \nu + n  e_{n+1}  \wedge e_{1} \wedge e_{2} \wedge \cdots \wedge e_{n},
\end{equation}
where  $\hat \alpha $ is the mean curvature of $M_t$ in  $\mathbb S^{n+1}_s(1)$.
\end{Lemma}

\proof 
Let $M_t$ be a pseudo--Riemannian hypersurface with index $t$ in \linebreak
$\mathbb S^{n+1}_s(1)\subset\mathbb E^{n+2}_s$.
Let $e_{1},  \ldots, e_{n+1}, e_{n+2}$ be a local orthonormal frame field on $M_t$ in $\mathbb E^{n+2}_s$ 
such that $e_{1}, e_{2}, \ldots, e_{n}$ are tangent to $M_t$, and $e_{n+1}, e_{n+2}= {\bf x}$ are normal to $M_t$, where ${\bf x}$ is the position vector of $M_t$. 
Since $e_{n+2} = {\bf x}$ is parallel in the normal bundle of $M_t$ in $\mathbb E^{n+2}_s$ and the codimension of $M_t$ in $\mathbb E^{n+2}_s$ is two, $e_{n+1}$ is   parallel too, i.e., $D  e_{n+1}=0$.
We put $\bar \nu = e_{n+1}  \wedge e_{1} \wedge e_{2} \wedge \cdots \wedge e_{n}$. 
Now we will compute $\Delta \bar \nu$. By differentiating $\bar \nu$ we get 
%%%%%%%%%%%%%%%%%
\begin{equation}\label{differen-1}
e_i \bar  \nu =  -  \varepsilon_i e_{n+1} \wedge  e_{1} \wedge \cdots \wedge e_{i-1} 
\wedge  {\bf x} \wedge e_{i+1} \wedge \cdots \wedge e_{n}.
 \end{equation}
Since   ${\nabla}_{e_i}e_i= \sum_{j=1}^{n} \varepsilon_j \omega_{ij}(e_i) e_j$, we have 
\begin{align} \label{diff-con}
({\nabla}_{e_i}e_i) \bar  \nu =- \sum_{j=1}^n  \omega_{ij}(e_i)  \, e_{n+1} \wedge e_{1} \wedge \cdots  
\wedge e_{j-1} \wedge  {\bf x} \wedge e_{j+1} \wedge   \cdots \wedge e_{n}.
\end{align}  
Differentiating $e_i \bar \nu $ in  \eqref{differen-1} we obtain that
%%%%
\begin{align} \label{diffent-2}
\begin{split}
e_ie_i \bar  \nu = - \varepsilon_i\bar \nu - h_{ii}^{n+1} \tilde \nu + \sum_{j=1}^n  \omega_{ji}(e_i)  
 e_{n+1} \wedge e_{1} \wedge  \cdots \wedge
\underbrace{{\bf x}}_{j-th} \wedge \cdots \wedge e_{n}.
\end{split}
\end{align}  
Considering $n \hat \alpha = \sum_{i=1}^n\varepsilon_i h^{n+1}_{ii}$  we have 
%%%%%%%%%%%%%%%%
\begin{align}\label{laplce-nu-bar}
\begin{split}
\Delta \bar  \nu = &  \sum_{i=1}^n  \varepsilon_i(\nabla_{e_i}e_i -e_i e_i ) \bar \nu   \\
                  = & n \hat \alpha   \tilde \nu    +n \bar \nu - 
					\sum_{i,j=1}^n \varepsilon_i (\omega_{ij}(e_i) + \omega_{ji}(e_i))   
					\, e_{n+1} \wedge e_{1} \wedge \cdots \wedge 
\underbrace{{\bf x}}_{j-th} \wedge \cdots \wedge e_{n}              
\end{split}
\end{align}
which gives \eqref{laplce-lem-2} as $\omega_{ij} + \omega_{ji} = 0$. 
\eproof

\theoremstyle{definiton}
\begin{example}\cite{Lucas}(Totally umbilical hypersurfaces in $\mathbb{S}^{n+1}_s$) \\
\label{extotumb} 
Totally umbilical hypersurfaces in $\mathbb{S}^{n+1}_s$ are obtained 
as the intersection of $\mathbb{S}^{n+1}_s$ with a hyperplane of $\mathbb{E}^{n+2}_q$, 
and the causal character of the hyperplane determines the type of the hypersurface. 
Let $a\in\mathbb{E}^{n+2}_q$ be a nonzero constant vector with $\langle a, a\rangle\in\{-1,0,1\}$.
For every $\tau\in\mathbb{R}$ with $\langle a, a\rangle-\tau^2\neq 0$, the set
\begin{equation}
\nonumber
M_t=\{{\bf x}\in\mathbb{S}^{n+1}_s:\langle {\bf x}, a\rangle=\tau \}
\end{equation}  
is a totally umbilical hypersurface with the unit normal vector in 
$\mathbb{S}^{n+1}_s\subset\mathbb{E}^{n+2}_s$
\begin{equation}
\nonumber
N=\frac{1}{\sqrt{|\langle a,a\rangle-\tau^2|}}(a-\tau{\bf x}),
\end{equation}
and the shape operator 
\begin{equation}
\label{shp}
A_N=\frac{\tau}{\sqrt{|\langle a,a\rangle-\tau^2|}}I_t,
\end{equation}
where $I_t$ is the identity transformation on $T_pM_t$ at $p\in M_t$.

Then, $M_t$ has constant curvature
\begin{equation}
\nonumber
K=1+\frac{\tau^2}{\langle a,a\rangle-\tau^2}.
\end{equation}
Hence, there are three different possibilities according to the causal character of the hyperplane 
as seen below:
\begin{itemize}
\item [(i).] If $\langle a, a\rangle =-1$, 
then $M_{t-1}=\mathbb{S}^n_{t-1}\left(\frac{1}{1+\tau^2}\right)\subset\mathbb{S}^{n+1}_{t}$ 
is a pseudo $n$--sphere of constant curvature $\frac{1}{1+\tau^2}$ for every $\tau$. 

\item [(ii).] Let $\langle a,a\rangle =1$. 
For $|\tau|<1$, $M_t=\mathbb{S}^{n}_t\left(\frac{1}{1-\tau^2}\right)\subset\mathbb{S}^{n+1}_t$ 
is a pseudo $n$--sphere of constant curvature $\frac{1}{1-\tau^2}$. 
For $|\tau|>1$, $M_{t-1}=\mathbb{H}^{n}_{t-1}\left(-\frac{1}{1-\tau^2}\right)\subset\mathbb{S}^{n+1}_{t}$ 
is a pseudo--hyperbolic $n$--space of constant curvature $\frac{1}{\tau^2-1}$. 

\item [(iii).] If $\langle a, a\rangle=0$, then $M_t$ is a flat totally umbilical hypersurface 
in $\mathbb{S}^{n+1}_{t+1}$ which is called the pseudo--horosphere in $\mathbb{S}^{n+1}_{t+1}$.
It follows from \eqref{shp} that $M_t$ has constant mean curvature $\hat{\alpha}$
with $|\hat{\alpha}|=1$.
\end{itemize}  
\end{example}

It is seen from the above example that there are totally umbilical hypersurfaces in 
$\mathbb{S}^{n+1}_s$ which are flat or non--totally geodesic.

\begin{Theorem}
\label{1typenonmass}
An n--dimensional pseudo--Riemannian submanifold $M_t$ with index $t$ and 
non--null mean curvature vector of a pseudo--sphere $\mathbb{S}^{m-1}_s(1)\subset\mathbb{E}^m_s$ 
has 1--type pseudo-spherical Gauss map 
with a nonzero constant component in its spectral decomposition 
if and only if $M_t$ is an open part of a non--flat, 
non--totally geodesic and totally umbilical pseudo--Riemannian hypersurface 
of a totally geodesic pseudo--sphere 
$\mathbb{S}^{n+1}_{s^*}(1)\subset\mathbb{S}^{m-1}_s(1)$, ($s^*=t\leq s$ or $s^*=t+1\leq s$), 
that is, it is an open portion of $\mathbb{S}^n_t(c)\subset\mathbb{S}^{n+1}_t(1)$ 
of curvature $c$ for $c>1$ or 
$\mathbb{S}^n_t(c)\subset\mathbb{S}^{n+1}_{t+1}(1)$ of curvature $c$ for $0<c<1$ 
or $\mathbb{H}^n_t(-c)\subset\mathbb{S}^{n+1}_{t+1}(1)$ of curvature $-c$ for $c>0$. 
\end{Theorem}

\proof
Let ${\bf x}: M_t\longrightarrow\mathbb{S}^{m-1}_s(1)\subset\mathbb{E}^m_s$ be an isometric immersion 
of an $n$--dimensional pseudo--Riemannian manifold $M_t$ with non--null mean curvature vector 
$\hat{H}$ in $\mathbb{S}^{m-1}_s(1)$ into 
$\mathbb{S}^{m-1}_s(1)\subset\mathbb{E}^m_s$. 
If the pseudo--spherical Gauss map $\tilde\nu$ is of 1--type with nonzero constant component in its spectral decomposition, 
then $\Delta\tilde\nu=\lambda_p(\tilde\nu-c)$ for some nonzero constant vector $c$ and a nonzero real number 
$\lambda_p$
Thus we have
\begin{equation}\label{non-mass-condition}
(\Delta \tilde \nu)_i = \lambda_p (\tilde{\nu})_i, 
\end{equation}
where $(\cdot)_i = e_i(\cdot)$.
On the other hand, by a direct long calculation, we get from \eqref{calc-laplce-lem} that 
\begin{align}
\label{deriv-of-lap} 
\nonumber  
e_i(\Delta &\tilde{\nu})
=(\|\hat{h}\|^2)_i\tilde{\nu}\\
&+\|\hat{h}\|^2\sum_{r=n+1}^{m-1} \sum_{k=1}^n \varepsilon_r h_{ik}^r {\bf x}\wedge e_1\wedge\cdots\wedge e_{k-1} \wedge e_r 
\wedge e_{k+1}\wedge\cdots\wedge e_n\notag\\ 
&+2n D_{e_i}\hat{H}\wedge e_1\wedge\cdots\wedge e_n
 +n\sum_{k=1}^n \sum_{r=n+1}^{m-1} \varepsilon_r {h}_{ik}^r \hat{H}\wedge e_1\wedge \cdots\wedge 
\underbrace{e_r}_{k-th}\wedge\cdots\wedge e_n\notag\\
&-n \sum_{k=1}^n \varepsilon_i\delta_{ik}\hat{H}\wedge e_1\wedge \cdots\wedge e_{k-1} 
     \wedge {\bf x} \wedge e_{k+1}\wedge\cdots\wedge e_n\notag\\
&-n\sum_{j,k,\ell=1 \atop j\neq k}^n \varepsilon_\ell\omega_{j\ell}(e_i) {\bf x} \wedge e_1\wedge\cdots\wedge 
    \underbrace{e_\ell}_{j-th}\wedge\cdots\wedge \underbrace{D_{e_k}\hat{H}}_{k-th}\wedge\cdots\wedge e_n\notag\\
&-n\sum_{j,k=1\atop j\neq k }^n \; \sum_{r=n+1}^{m-1} \varepsilon_r {h}_{ij}^r {\bf x}\wedge e_1 
      \wedge\cdots\wedge \underbrace{e_r}_{j-th} \wedge \cdots\wedge 
			 \underbrace{D_{e_k}\hat{H}}_{k-th}\wedge\cdots\wedge e_n\notag\\
&+n\sum_{k=1}^n \varepsilon_k \langle A_{D_{e_k}\hat{H}}(e_i),e_k \rangle {\bf x}\wedge e_1\wedge\cdots\wedge e_n \notag\\
&- n \sum_{k=1}^n {\bf x} \wedge e_1\wedge\cdots\wedge \underbrace{D_{e_i}D_{e_k}\hat{H}}_{k-th}\wedge \cdots\wedge e_n \notag\\
&+\sum_{r,s=n+1\atop r<s}^{m-1} \sum_{j,k=1\atop j\neq k}^n \varepsilon_r\varepsilon_s\{e_i(R_{sjk}^r) {\bf x}+ R_{sjk}^r {e_i}\} 
      \wedge e_1\wedge\cdots\wedge \underbrace{e_r}_{j-th}\wedge\cdots\wedge
\underbrace{e_s}_{k-th}\wedge\cdots\wedge e_n\notag\\
&+\sum_{r,s=n+1\atop r<s}^{m-1} \sum_{j,k,\ell=1 \atop j,k,\ell\neq}^n \varepsilon_r\varepsilon_s R_{sjk}^r  
       \Big  \{ \sum_{h=1}^n \varepsilon_h\omega_{\ell h}(e_i) {\bf x}\wedge e_1\wedge\cdots\wedge 
			  \underbrace{e_h}_{\ell-th}\wedge\cdots\wedge\underbrace{e_r}_{j-th}\wedge\cdots\wedge \underbrace{e_s}_{k-th}      \notag \\
&\hskip 5mm \wedge\cdots\wedge e_n +\sum_{t=n+1}^{m-1} \varepsilon_t{h}_{i\ell}^t {\bf x}\wedge e_1\wedge\cdots\wedge 
      \underbrace{e_t}_{\ell-th}\wedge\cdots\wedge\underbrace{e_r}_{j-th}\wedge\cdots    \wedge  \underbrace{e_s}_{k-th}\wedge\cdots\wedge e_n \Big \} \notag\\      
& + \sum_{r,s=n+1}^{m-1} \sum_{j,k,\ell=1\atop j\neq k}^n \varepsilon_\ell\varepsilon_r\varepsilon_s R_{sjk}^r {h}_{i\ell}^s {\bf x} 
      \wedge e_1\wedge\cdots\wedge \underbrace{e_\ell}_{j-th}\wedge\cdots\wedge
        \underbrace{e_r}_{k-th}\wedge\cdots\wedge e_n\notag \\           
&-\sum_{r,s,t=n+1}^{m-1} \sum_{j,k=1\atop j\neq k}^n \varepsilon_r\varepsilon_s\varepsilon_t  R_{sjk}^r\omega_{st}(e_i) {\bf x}\wedge e_1\wedge\cdots\wedge \underbrace{e_t}_{j-th} 
\wedge\cdots\wedge\underbrace{e_r}_{k-th} \wedge\cdots\wedge e_n.  
\end{align}

\textit{Case (a)}: $\hat{H}=0$. So, equation \eqref{deriv-of-lap} reduces to
\begin{align} 
\label{deriv-of-lap-H=0}
e_i(\Delta &\tilde{\nu})
=(\|\hat{h}\|^2)_i\tilde{\nu}+ \|\hat{h}\|^2\sum_{r=n+1}^{m-1} 
     \sum_{k=1}^n \varepsilon_r h_{ik}^r {\bf x}\wedge e_1\wedge\cdots\wedge 
     \underbrace{e_r}_{k-th}\wedge\cdots\wedge e_n \notag\\
&+\sum_{r,s=n+1\atop r<s}^{m-1} \sum_{j,k=1\atop j\neq k}^n \varepsilon_r\varepsilon_s\{e_i(R_{sjk}^r) {\bf x}+ 
      R_{sjk}^r {e_i}\}\wedge e_1\wedge\cdots\wedge \underbrace{e_r}_{j-th}\wedge\cdots\wedge
      \underbrace{e_s}_{k-th}\wedge\cdots\wedge e_n\notag\\
&+\sum_{r,s=n+1\atop r<s}^{m-1} \sum_{j,k,\ell\atop j,k,\ell\neq}^n \varepsilon_r\varepsilon_s R_{sjk}^r  
       \Big  \{ \sum_{h=1}^n \varepsilon_h\omega_{\ell h}(e_i) {\bf x}\wedge e_1\wedge\cdots\wedge 
			  \underbrace{e_h}_{\ell-th}\wedge\cdots\wedge\underbrace{e_r}_{j-th}\wedge\cdots\wedge \underbrace{e_s}_{k-th}      \notag \\
&\hskip 5mm \wedge\cdots\wedge e_n +\sum_{t=n+1}^{m-1} \varepsilon_t{h}_{i\ell}^t {\bf x}\wedge e_1\wedge\cdots\wedge 
      \underbrace{e_t}_{\ell-th}\wedge\cdots\wedge\underbrace{e_r}_{j-th}\wedge\cdots    \wedge  \underbrace{e_s}_{k-th}\wedge\cdots\wedge e_n \Big \} \notag     \\
&+ \sum_{r,s=n+1}^{m-1} \sum_{j,k,\ell=1\atop j\neq k}^n \varepsilon_\ell\varepsilon_r\varepsilon_s R_{sjk}^r {h}_{i\ell}^s {\bf x} 
       \wedge e_1\wedge\cdots\wedge \underbrace{e_\ell}_{j-th}\wedge\cdots\wedge
       \underbrace{e_r}_{k-th}\wedge\cdots\wedge e_n\notag\\
&-\sum_{r,s,t=n+1}^{m-1} \sum_{j,k=1\atop j\neq k}^n \varepsilon_t\varepsilon_r\varepsilon_s R_{sjk}^r\omega_{st}(e_i) {\bf x}\wedge e_1\wedge\cdots\wedge \underbrace{e_t}_{j-th} 
\wedge\cdots\wedge\underbrace{e_r}_{k-th}\wedge\cdots\wedge e_n. 
\end{align}
Comparing \eqref{dif1}, \eqref{non-mass-condition} and \eqref{deriv-of-lap-H=0} 
we can see that $(\|\hat{h}\|^2)_i =R_{sjk}^r = 0$ for all $i,j,k=1,2,\dots,n$ and $r,s=n+1,\dots,m-1$. 
Thus, according to Theorem \ref{1typesubmanifold} the pseudo--spherical Gauss map $\tilde\nu$ 
of a pseudo--Riemannian submanifold $M_t$ is of 1--type with $c=0$ which is a contradiction.  

\textit{Case (b)}: $\hat{H}\neq 0$. Since the term 
$D_{e_i}\hat{H}\wedge e_1\wedge\cdots\wedge e_n$  
appears only in $e_i(\Delta \tilde{\nu})$, not in $e_i(\tilde \nu)$, we obtain from 
\eqref{dif1}, \eqref{non-mass-condition} and \eqref{deriv-of-lap} that $D{\hat H}=0$ 
which implies that $M_t$ has parallel mean curvature vector $\hat{H}$ in $\mathbb S^{m-1}_s(1)$, 
and hence, $\langle\hat{H},\hat{H}\rangle$ is a nonzero constant. 
Thus, equation \eqref{deriv-of-lap} reduces to
\begin{align} 
\label{deriv-of-lap-H-no-0}
e_i(\Delta &\tilde{\nu})
=(\|\hat{h}\|^2)_i\tilde{\nu} + 
   \|\hat{h}\|^2\sum_{r=n+1}^{m-1} \sum_{k=1}^n \varepsilon_r h_{ik}^r {\bf x} 
		\wedge e_1\wedge\cdots\wedge\underbrace{e_r}_{k-th}\wedge\cdots\wedge e_n\notag\\
&+ n\sum_{k=1}^n \sum_{r=n+1}^{m-1} \varepsilon_r {h}_{ik}^r \hat{H}\wedge e_1\wedge 
     \cdots\wedge e_{k-1} \wedge e_r \wedge  e_{k+1}\wedge\cdots\wedge e_n\notag\\%\displaybreak
&-n \sum_{k=1}^n \varepsilon_i\delta_{ik}\hat{H}\wedge e_1\wedge \cdots\wedge e_{k-1} 
     \wedge {\bf x} \wedge  e_{k+1}\wedge\cdots\wedge e_n\notag\\ 
&+\sum_{r,s=n+1\atop r<s}^{m-1} \sum_{j,k=1\atop j\neq k}^n \varepsilon_r\varepsilon_s\{e_i(R_{sjk}^r) {\bf x}+ 
     R_{sjk}^r {e_i}\}\wedge e_1\wedge\cdots\wedge \underbrace{e_r}_{j-th}\wedge\cdots\wedge
      \underbrace{e_s}_{k-th}\wedge\cdots\wedge e_n\notag\\
&+\sum_{r,s=n+1\atop r<s}^{m-1} \sum_{j,k,\ell\atop j,k,\ell\neq}^n \varepsilon_r\varepsilon_s R_{sjk}^r  
       \Big  \{ \sum_{h=1}^n \varepsilon_h\omega_{\ell h}(e_i) {\bf x}\wedge e_1\wedge\cdots\wedge 
			  \underbrace{e_h}_{\ell-th}\wedge\cdots\wedge\underbrace{e_r}_{j-th}\wedge\cdots\wedge \underbrace{e_s}_{k-th}      \notag \\
&\hskip 5mm \wedge\cdots\wedge e_n +\sum_{t=n+1}^{m-1} \varepsilon_t{h}_{i\ell}^t {\bf x}\wedge e_1\wedge\cdots\wedge 
      \underbrace{e_t}_{\ell-th}\wedge\cdots\wedge\underbrace{e_r}_{j-th}\wedge\cdots    \wedge  \underbrace{e_s}_{k-th}\wedge\cdots\wedge e_n \Big \} \notag     \\
 & + \sum_{r,s=n+1}^{m-1} \sum_{j,k,\ell=1\atop j\neq k}^n \varepsilon_\ell\varepsilon_r\varepsilon_sR_{sjk}^r {h}_{i\ell}^s {\bf x} 
 \wedge e_1\wedge\cdots\wedge \underbrace{e_\ell}_{j-th}\wedge\cdots\wedge
 \underbrace{e_r}_{k-th}\wedge\cdots\wedge e_n\notag\\
&-\sum_{r,s,t=n+1}^{m-1} \sum_{j,k=1\atop j\neq k}^n \varepsilon_r\varepsilon_s\varepsilon_t R_{sjk}^r\omega_{st}(e_i){\bf x}\wedge e_1\wedge\cdots\wedge \underbrace{e_t}_{j-th} \wedge\cdots\wedge\underbrace{e_r}_{k-th} \wedge\cdots\wedge e_n. 
\end{align}
From  \eqref{dif1}, \eqref{non-mass-condition} and \eqref{deriv-of-lap-H-no-0} 
we obtain that $\|\hat h\|^2$ is constant which means that the scalar curvature of $M_t$ is constant. 
Also, we have 
\begin{align} 
\label{Equ-1}
 &   \|\hat{h}\|^2\sum_{r=n+1}^{m-1} \sum_{k=1}^n \varepsilon_r h_{ik}^r {\bf x} 
       \wedge e_1\wedge\cdots\wedge   e_{k-1} \wedge e_r \wedge  e_{k+1} \wedge\cdots\wedge e_n\notag\\
&-n \sum_{k=1}^n \varepsilon_i\delta_{ik}\hat{H}\wedge e_1\wedge \cdots\wedge e_{k-1} 
     \wedge {\bf x} \wedge  e_{k+1}\wedge\cdots\wedge e_n\notag\\
&+ \sum_{r,s=n+1}^{m-1} \sum_{j,k,\ell=1\atop j\neq k}^n \varepsilon_\ell\varepsilon_r\varepsilon_s R_{sjk}^r {h}_{i\ell}^s {\bf x} 
       \wedge e_1\wedge\cdots\wedge \underbrace{e_\ell}_{j-th}\wedge\cdots\wedge
        \underbrace{e_r}_{k-th}\wedge\cdots\wedge e_n\notag\\
&= \lambda_p \sum_{r=n+1}^{m-1} \sum_{k=1}^n \varepsilon_r h_{ik}^r {\bf x}\wedge e_1\wedge\cdots\wedge e_{k-1} 
      \wedge e_r \wedge  e_{k+1}\wedge\cdots\wedge e_n
\end{align}
and 
\begin{align} 
\label{Equ-2}
& n\sum_{k=1}^n \sum_{r=n+1}^{m-1} \varepsilon_r {h}_{ik}^r \hat{H}\wedge e_1 
 \wedge \cdots\wedge  e_{k-1} \wedge {e_r} \wedge e_{k+1} 
 \wedge\cdots\wedge e_n\notag\\
&+\sum_{r,s=n+1 \atop r<s}^{m-1} \sum_{j,k=1\atop j\neq k}^n  \varepsilon_r\varepsilon_s R_{sjk}^r {e_i}\wedge e_1\wedge\cdots\wedge 
       \underbrace{e_r}_{j-th}\wedge\cdots\wedge
        \underbrace{e_s}_{k-th}\wedge\cdots\wedge e_n  = 0. 
\end{align}
We put $\hat{H}=\varepsilon_{n+1}\hat{\alpha}e_{n+1}$. It follows from \eqref{Equ-2} that 
$R^r_{sjk} = 0$ for $r, s \geq n + 2$ and $j, k = 1,\ldots n$. 
Also, we  find $R^{n+1}_{sjk} = 0$ from $D\hat{H} = 0$. 
Thus, the normal connection of $M_t$ in $\mathbb S^{m-1}_s(1)$ is flat.
Therefore, \eqref{Equ-2} yields 
\begin{align} 
\label{Equ-3}
n \varepsilon_{n+1} \hat \alpha  \sum_{k=1}^n \sum_{r=n+1}^{m-1} \varepsilon_r {h}_{ik}^r e_{n+1}
\wedge e_1\wedge \cdots\wedge e_{k-1} \wedge  e_r \wedge e_{k+1} \wedge \cdots\wedge e_n =0
\end{align}
which gives $h_{ik}^r=0$ for $r\geq n+2$ and $i, k=1,2,\dots, n$, and hence the first normal space $\mbox{Im}h$ is spanned by $e_{n+1}$. Therefore, by the reduction theorem of Erbarcher, 
we conclude that $M_t$ is contained in a totally geodesic pseudo--sphere $\mathbb{S}^{n+1}_{s^*}(1)$ 
in $\mathbb{S}^{m-1}_s(1)\subset\mathbb{E}^m_s$ 
for $s^*=t$ or $s^*=t+1$. Also, from \eqref{Equ-1} we have
\begin{equation}
\label{pseudo-umbilical-eq}
\varepsilon_i (\|\hat{h}\|^2  -  \lambda_p )   h_{ik}^{n+1}+ n\delta_{ik}{\hat \alpha } =0 
\end{equation} 
for $i, k = 1,2,\ldots n$. If $\lambda_p = \|\hat{h}\|^2  $, then \eqref{pseudo-umbilical-eq} gives 
$ \hat \alpha = 0$ which is a contradiction. So, $\lambda_p \not  = \|\hat{h}\|^2$, and 
by taking the sum of \eqref{pseudo-umbilical-eq} for  $i= k$ and $i$ from 1 up to $n$ we get 
$n\hat \alpha (n +  \|\hat{h}\|^2  -  \lambda_p )=0$ which gives  $\lambda_p= n +  \|\hat{h}\|^2$. 
Thus, we get from \eqref{pseudo-umbilical-eq} that 
$h_{ii}^{n+1} = \varepsilon_i\hat \alpha \not = 0$ and $h_{ik}^{n+1}=0$ for $i\neq k$. 
Therefore, $M_t$ is a non--totally geodesic and totally umbilical pseudo--Riemannian hypersurface of 
$\mathbb{S}^{n+1}_{s^*}(1)\subset\mathbb{S}^{m-1}_s(1)$. 
From \eqref{scalar-curv-sphere} we get $S=n(n-1)(1+\varepsilon_{n+1}\hat\alpha^2)\neq 0$ 
because of $\lambda_p=n(1+\varepsilon_{n+1}\hat\alpha^2)\neq 0$. 
Consequently, it follows from Example \ref{extotumb} that $M_t$ is an open part 
of $\mathbb{S}^n_t(c)\subset\mathbb{S}^{n+1}_{t+1}(1)$ of curvature $c$ for $0<c<1$ or 
$\mathbb{S}^n_t(c)\subset\mathbb{S}^{n+1}_t(1)$ of curvature $c$ for $c>1$ or 
$\mathbb{H}^n_t(-c)\subset\mathbb{S}^{n+1}_{t+1}(1)$ of curvature $-c$ for $c>0$.

Conversely, let $M_t$ be an open part of a nonflat, non--totally geodesic and totally umbilical 
hypersurface of a totally geodesic pseudo--sphere 
$\mathbb{S}^{n+1}_{s^*}(1)\subset\mathbb{S}^{m-1}_s$. 
Without loss of generality, we assume that $M_t$ is immersed in 
$\mathbb{S}^{n+1}_{s^*}(1)\subset\mathbb{E}^{n+2}_{s^*}$. 
Let $e_1, \dots, e_{n+1}, e_{n+2}={\bf x}$ be a local orthonormal frame 
on $M_t$ in $\mathbb{E}^{n+2}_s$ 
such that $e_1,\dots, e_n$ are tangent to $M_t$, and $e_{n+1},e_{n+2}$ are normal to $M_t$, 
where ${\bf x}$ is the position vector of $M_t$. 
Since $M_t$ is a hypersurface of $\mathbb{S}^{n+1}_{s^*}(1)$, 
the normal bundle of $M_t$ in $\mathbb{E}^{n+2}_{s^*}$ is flat. 
Also, the mean curvature vector $\hat{H}=\varepsilon_{n+1}\hat\alpha e_{n+1}$ is parallel 
in $\mathbb{E}^{n+2}_{s^*}$.
As $M_t$ is totally umbilical, we get 
$\|\hat{h}\|^2=\varepsilon_{n+1}n\hat{\alpha}^2$ and hence, 
from \eqref{calc-laplce-lem} we have
\begin{equation}
\label{first-lap-hyper-sur}
\Delta  \tilde \nu =  n\varepsilon_{n+1}\hat\alpha(\hat\alpha\tilde\nu
+e_{n+1}\wedge e_1\wedge e_2\wedge\cdots\wedge e_n).
\end{equation}
Now, if we put 
$$
c = \frac{1}{1+\varepsilon_{n+1}\hat \alpha^2}(\tilde \nu - \varepsilon_{n+1}\hat \alpha e_{n+1}  \wedge e_{1} \wedge e_{2} \wedge \cdots \wedge e_{n}),
$$
$$
\tilde \nu_p = \frac{\varepsilon_{n+1}\hat \alpha }{1+\varepsilon_{n+1}\hat \alpha^2}(\hat \alpha \tilde \nu  +  
e_{n+1}  \wedge e_{1} \wedge e_{2} \wedge \cdots \wedge e_{n}),
$$
where ${1+\varepsilon_{n+1}\hat \alpha^2}\neq 0$ because $M_t$ is nonflat hypersurface 
in $\mathbb{S}^{n+1}_{s^*}(1)$
(note that for a flat totally umbilical hypersurface in $\mathbb{S}^{n+1}_{s^*}(1)$, 
$\varepsilon_{n+1}=-1$ and $\hat{\alpha}^2=1$), then we have $\tilde\nu=c+\tilde\nu_p$. 
As $\hat \alpha$ is constant it is easily seen that 
$e_i(c) = 0 $, $i= 1, \dots, n$, i.e.,  $c$ is a constant vector.
Using \eqref{laplce-lem-2} and \eqref{first-lap-hyper-sur},  
we obtain from  a direct computation that   
$\Delta \tilde \nu_p = n(1+\varepsilon_{n+1}\hat \alpha^2) \tilde \nu_p$. 
Therefore, the pseudo--spherical Gauss map $\tilde \nu$ is of 1--type with nonzero constant component 
in its spectral decomposition.
\eproof

\begin{Definition}
A pseudo--Riemannian submanifold is said to have biharmonic pseudo--spherical Gauss map 
$\tilde\nu$ if it satisfies $\Delta\tilde\nu\neq 0$ and 
$\Delta^2\tilde\nu=0$ identically.
\end{Definition}

\begin{Theorem}
A pseudo--horosphere $M_t$ in the pseudo--sphere $\mathbb{S}^{n+1}_s(1)\subset\mathbb{E}^{n+2}_s$ has biharmonic pseudo--spherical Gauss map.
\end{Theorem}

\proof
Suppose that $M_t$ is a pseudo--horosphere with index $t$ in \linebreak
$\mathbb{S}^{n+1}_s(1)\subset\mathbb{E}^{n+2}_s$. 
As seen in Example \ref{extotumb}, 
for given any null vector $a\in\mathbb{E}^{n+2}_s$ and a real number $\tau\neq 0$, 
an $n$--dimensional pseudo--horosphere in $\mathbb{S}^{n+1}_s(1)$ is defined by
\begin{equation}
\nonumber
M_t=\{{\bf x}\in\mathbb{S}^{n+1}_s(1):\langle {\bf x},a \rangle=\tau\}.
\end{equation}
The unit normal vector $e_{n+1}$ of $M_t$ in $\mathbb{S}^{n+1}_s(1)$ is chosen as
$e_{n+1}=\frac{1}{\tau}(a-\tau{\bf x})$ with $\varepsilon_{n+1}=-1$. 
Thus, $A_{n+1}= I_t$ from which we get the mean curvature 
$\hat{\alpha}= 1$ and $\|\hat{h}\|^2=-n$. Thus, \eqref{calc-laplce-lem} becomes
\begin{equation}
\label{laphoro}
\Delta\tilde\nu=-n\tilde\nu-n e_{n+1}\wedge e_1\wedge e_2\wedge\cdots\wedge e_n\neq 0. 
\end{equation}
On the other hand, from \eqref{laplce-lem-2} we have 
\begin{equation}
\label{laphoro1}
\Delta({e_{n+1}\wedge e_1\wedge e_2\cdots\wedge e_n})=n\tilde\nu
+n e_{n+1}\wedge e_1\wedge e_2\wedge\cdots\wedge e_n=-\Delta\tilde\nu.
\end{equation}
Hence, we obtain that
\begin{equation}
\label{laphoro2}
\Delta^2\tilde\nu=-n\Delta\tilde\nu-n\Delta(e_{n+1}\wedge e_1\wedge e_2\wedge\cdots\wedge e_n)
=-n\Delta\tilde\nu+n\Delta\tilde\nu=0
\end{equation}
which implies that the pseudo--spherical Gauss map $\tilde\nu$ of the pseudo--horosphere is biharmonic. 
\eproof

\begin{Definition}
A surface of a pseudo--sphere is marginally trapped 
if and only if its mean curvature vector is lightlike at each point on the surface.
\end{Definition}

By Theorem 6.1 in \cite{chenveken2}, B.-Y. Chen and J. Van der Veken classified the marginally trapped surfaces 
with parallel mean curvature vector in $\mathbb{S}^4_1(1)\subset\mathbb{E}^4_1$, 
and they obtained eight types of surfaces.   
It follows from the proof of this theorem the surface given in the next example is 
the only totally umbilical marginally trapped surface.

\begin{example} 
%\cite{chenveken2}
\label{surfnonmassquant}
Consider the following spacelike surface $M$ in $\mathbb{S}^4_1(1)\subset\mathbb{E}^5_1$ defined by
\begin{equation}
\label{surfnonmass}
{\bf x}(u, v)=(1, \sin u, \cos u\cos v, \cos u\sin v, 1).
\end{equation}
We choose an orthonormal moving frame field $e_1, e_2, \dots, e_5$ on $M$ 
such that $e_1,e_2$ are tangent to $M$ and $e_3,e_4,e_5$ are normal to $M$
in $\mathbb{E}^5_1$ as follows:
\begin{eqnarray*}
\label{tangent12} &e_1=\frac{\partial}{\partial u}, \quad  
e_2=\frac{1}{\cos u}\frac{\partial}{\partial v},\;\; \cos u\neq 0,\\
\label{normal3} &e_3=\frac{1}{\sqrt{2}}(-1, \sin u, \cos u\cos v, \cos u\sin v, 0),\\
\label{normal45} &e_4=\frac{1}{\sqrt{2}}(1, \sin u, \cos u\cos v, \cos u\sin v, 2), \quad e_5={\bf x},
\end{eqnarray*}
where $\varepsilon_1=\varepsilon_2=\varepsilon_3=-\varepsilon_4=\varepsilon_5=1$. 
By a direct calculation, we have the shape operators and the connection forms as
\begin{eqnarray*}
A_3=A_4=-\frac{1}{\sqrt{2}}\mbox{I},\;\; \omega_{12}(e_1)=0,  \;\; \omega_{12}(e_2)=-\tan u,\;\; \omega_{34}=0,
\end{eqnarray*} 
where I is the $2\times 2$ identity matrix. Thus, the mean curvature vector in $\mathbb{S}^4_1$ is given by 
$\hat{H}=-\frac{1}{\sqrt{2}}(e_3-e_4)$ which is parallel and null, and the normal curvature $K^D\equiv 0$ and 
the Gauss curvature $K=1$.
Moreover, from \eqref{laplce-nu-bar} the Laplacian of the pseudo-spherical Gauss map of $M$ is 
$\Delta\tilde\nu=2\hat{H}\wedge e_1\wedge e_2\neq 0$. 
Then, $\tilde\nu$ has the spectral decomposition $\tilde\nu=c+\tilde\nu_p$ with 
$c=\tilde\nu+\frac{1}{\sqrt{2}}(e_3\wedge e_1\wedge e_2 - e_4\wedge e_1\wedge e_2)$, 
$\tilde\nu_p=-\frac{1}{\sqrt{2}}(e_3\wedge e_1\wedge e_2 - e_4\wedge e_1\wedge e_2)$ and 
$\Delta\tilde\nu_p=2\tilde\nu_p$. It can be seen easily that $c$ is a constant vector.  
Therefore, the pseudo--spherical Gauss map $\tilde\nu$ of $M$ is of 1--type with nonzero constant component in its spectral decomposition.
\end{example}

\begin{Theorem}
\label{1typenonmassS4_1}
A marginally trapped spacelike surface $M$ in the de Sitter space $\mathbb{S}^4_1(1)\subset\mathbb{E}^5_1$ 
has 1--type pseudo--spherical Gauss map with nonzero constant component in its spectral decomposition 
if and only if $M$ is congruent to an open part of the surface of curvature one given by \eqref{surfnonmass}.     
\end{Theorem}

\proof
Let ${\bf x}: M\longrightarrow\mathbb{S}^{4}_1(1)\subset\mathbb{E}^5_1$ be an isometric immersion of a marginally trapped spacelike surface $M$ into $\mathbb{S}^{4}_1(1)$. 
If the pseudo--spherical Gauss map $\tilde\nu$ is of 1--type with nonzero constant component
in its spectral decomposition, 
then we have $\Delta\tilde\nu=\lambda_p(\tilde\nu-c)$ for some nonzero constant vector $c$ and 
a real number $\lambda_p\neq 0$.  
Thus, we have
\begin{equation}\label{1non-mass-condition}
(\Delta \tilde \nu)_i = \lambda_p (\tilde{\nu})_i, 
\end{equation}
where $(\cdot)_i = e_i(\cdot)$.
On the other hand, by a direct calculation, we get from \eqref{calc-laplce-lem} that
\begin{align}
\label{1deriv-of-lap}
\nonumber  
e_i(\Delta &\tilde{\nu})
=\left((\|\hat{h}\|^2)_i+2\langle A_{D_{e_1} \hat{H}}, e_1\rangle +2\langle A_{D_{e_2} \hat{H}}, e_2\rangle \right)\tilde{\nu} \notag\\
&+\|\hat{h}\|^2\left(\sum_{r=3}^{4} \varepsilon_r h_{i1}^r {\bf x}\wedge e_r\wedge e_{2} 
+\sum_{r=3}^{4} \varepsilon_r h_{i2}^r {\bf x}\wedge e_1\wedge e_r \right)\notag\\ 
&+4 D_{e_i}\hat{H}\wedge e_1\wedge e_2
+2\sum_{r=3}^{4} \varepsilon_r {h}_{i1}^r \hat{H}\wedge e_r\wedge e_2
+2\sum_{r=3}^{4} \varepsilon_r {h}_{i2}^r \hat{H}\wedge e_1\wedge e_r\notag\\
&-2 \delta_{i1}\hat{H}\wedge {\bf x}\wedge e_2
-2 \delta_{i2}\hat{H}\wedge e_1 \wedge {\bf x} \notag\\
&+2 \omega_{12}(e_i) {\bf x} \wedge D_{e_1}\hat{H} \wedge e_1
-2 \omega_{12}(e_i) {\bf x} \wedge e_2\wedge D_{e_2}\hat{H} \notag\\
&-2\left(\sum_{r=3}^{4} \varepsilon_r {h}_{i2}^r {\bf x}\wedge D_{e_1}\hat{H}\wedge e_r
+\sum_{r=3}^{4} \varepsilon_r {h}_{i1}^r {\bf x}\wedge e_r \wedge D_{e_2}\hat{H}\right) \notag\\
&- 2\left( {\bf x} \wedge D_{e_i}D_{e_1}\hat{H}\wedge e_2 
+{\bf x} \wedge e_1\wedge D_{e_i}D_{e_2}\hat{H}\right)\notag\\
&+2\varepsilon_3\varepsilon_4\left(e_i(R_{412}^3) {\bf x}+ R_{412}^3 {e_i}\right) 
      \wedge e_3\wedge e_4\notag\\
& -2 \varepsilon_3\varepsilon_4 R_{412}^3 \left(\sum_{j=1}^2 h_{ji}^3 {\bf x}\wedge e_j\wedge e_4 
+ \sum_{j=1}^{2}  h_{ji}^4 {\bf x} \wedge e_3\wedge e_j\right).
\end{align}
Since the term $D_{e_i}\hat{H}\wedge e_1\wedge e_2$ appears only in $e_i(\Delta \tilde{\nu})$, 
not in $e_i(\tilde \nu)$, we have $D{\hat H}=0$ which implies 
that $M$ has parallel nonzero mean curvature vector in $\mathbb S^4_1(1)$.
Thus, from Lemma 2.2 in \cite{chenveken1} the surface $M$ has flat normal bundle in $\mathbb{S}^4_1(1)$, i.e., $R^3_{412}=0$. 
In this case, equation \eqref{1deriv-of-lap} becomes
\begin{align} 
\label{1deriv-of-lap-H-no-0}
e_i(\Delta &\tilde{\nu})
=(\|\hat{h}\|^2)_i\tilde{\nu} 
+ \|\hat{h}\|^2\left(\sum_{r=3}^{4} \varepsilon_r h_{i1}^r {\bf x}\wedge e_r\wedge e_{2} 
+\sum_{r=3}^{4} \varepsilon_r h_{i2}^r {\bf x}\wedge e_1\wedge e_r \right)\notag\\ 
&-2 \delta_{i1}\hat{H}\wedge {\bf x}\wedge e_2
-2 \delta_{i2}\hat{H}\wedge e_1 \wedge {\bf x} \notag\\
&+2\sum_{r=3}^{4} \varepsilon_r {h}_{i1}^r \hat{H}\wedge e_r\wedge e_2
+2\sum_{r=3}^{4} \varepsilon_r {h}_{i2}^r \hat{H}\wedge e_1\wedge e_r.
\end{align}
Hence, from \eqref{dif1}, \eqref{1non-mass-condition} and \eqref{1deriv-of-lap-H-no-0}
we obtain $\|\hat h\|^2$ is constant which implies that scalar curvature is constant
by \eqref{scalar-curv-sphere},
\begin{align} 
\label{1Equ-2}  
\sum_{r=3}^{4} \varepsilon_r {h}_{i1}^r \hat{H}\wedge e_r \wedge e_2
+ \sum_{r=3}^{4} \varepsilon_r {h}_{i2}^r \hat{H}\wedge e_1 \wedge e_r=0 
\end{align} 
and 
\begin{align} 
\label{1Equ-1}
&\|\hat{h}\|^2\left(\sum_{r=3}^{4} \varepsilon_r h_{i1}^r {\bf x} \wedge e_r\wedge e_2
+\sum_{r=3}^{4} \varepsilon_r h_{i2}^r {\bf x} \wedge e_1\wedge e_r\right)\notag\\
&-2 \delta_{i1}\hat{H}\wedge {\bf x} \wedge e_2
-2 \delta_{i2}\hat{H}\wedge e_1\wedge {\bf x} \notag\\
&= \lambda_p \left(\sum_{r=3}^{4} \varepsilon_r h_{i1}^r {\bf x}\wedge e_r\wedge e_2
+\sum_{r=3}^{4} \varepsilon_r h_{i2}^r {\bf x}\wedge e_1\wedge e_r\right).
\end{align}
Since the mean curvature vector $\hat{H}$ is lightlike, 
we have $(\mbox{tr} A_3)^2=(\mbox{tr} A_4)^2$.
As $M$ is a spacelike surface with flat normal bundle in $\mathbb{S}^4_1(1)$, 
the shape operators are simultaneously diagonalizable. 
Then, there exists an orthonormal frame $\{e_1, e_2, e_3, e_4\}$ on $M$ 
such that the shape operators and the second fundamental forms are given by
\begin{eqnarray}
\label{shapeop}
A_3=\left(
\begin{array}{cc}
h_{11}^3 & 0\\
0 & h_{22}^3
\end{array}
\right), 
\qquad 
A_4=\left(
\begin{array}{cc}
h_{11}^4 & 0\\
0 & h_{22}^4
\end{array}
\right), 
\end{eqnarray}
\begin{eqnarray}
\nonumber
\label{secondfun}
h(e_1,e_1)=\varepsilon_3(h_{11}^3 e_3 - h_{11}^4 e_4),\;
h(e_1, e_2)=0,\;
h(e_2,e_2)=\varepsilon_3(h_{22}^3 e_3 - h_{22}^4 e_4).
\end{eqnarray}
Hence, from the Gauss equation, we obtain 
\begin{equation}
\label{Gausscur}
K=1+\varepsilon_3(h_{11}^3h_{22}^3-h_{11}^4h_{22}^4).
\end{equation}
In this case, \eqref{1Equ-2} yields $h_{i1}^4\mbox{tr}A_3-h_{i1}^3\mbox{tr}A_4=0$ and 
$h_{i2}^4\mbox{tr} A_3-h_{i2}^3\mbox{tr}A_4=0$ for $i=1,2$. 
Since the shape operators are diagonalizable and $\hat{H}\neq 0$, these system of equations hold identically.
Also, it follows from \eqref{1Equ-1} that
\begin{equation}
\label{equ-1t}
(\|\hat{h}\|^2-\lambda_p)h^r_{ij}+\delta_{ij}\mbox{tr}A_r=0 
\end{equation}
for all $i,j=1,2$ and $r=3,4$. If $\lambda_p=\|\hat{h}\|^2$, then \eqref{equ-1t} gives 
$\mbox{tr}A_3=\mbox{tr}A_4=0$ which is a contradiction. So, $\lambda_p\neq \|\hat{h}\|^2,$ 
and from equation \eqref{equ-1t} we also have $h_{11}^3=h_{22}^3$ and $h_{11}^4=h_{22}^4$. 
Therefore, \eqref{shapeop} and \eqref{secondfun} reduce to 
\begin{equation}
\label{shapeop1}
A_3=h_{11}^3\mbox{I}\;\;\; \mbox{and}\;\;\; 
A_4=h_{11}^4\mbox {I},
\end{equation}
\begin{equation}
\label{secondfun1}
h(e_1,e_1)=h(e_2,e_2)=\varepsilon_3(h_{11}^3 e_3 - h_{11}^4 e_4),\;\;
h(e_1, e_2)=0,
\end{equation}
where I is $2\times 2$ identity matrix. Therefore, the surface $M$ is non--totally geodesic
and totally umbilical in $\mathbb{S}^4_1(1)$.
On the other hand, using the fact that 
$(\mbox{tr}A_3)^2=(\mbox{tr}A_4)^2$ in \eqref{Gausscur}, we get $K=1$. 
Thus, following the proof of Theorem 6.1 in \cite{chenveken2} it is seen that $M$ is congruent to an open part of the surface of curvature one whose position vector is given by \eqref{surfnonmass}.

The converse follows from Example \ref{surfnonmassquant}. 
\eproof

\textit{Acknowledgement.}
This work is a part of the first author's doctoral thesis which is 
partially supported by Istanbul Technical University.

\end{document}